\documentclass[a4paper,11pt]{article}
\usepackage{indentfirst, latexsym, bm}
\usepackage[numbers,sort&compress]{natbib}
\usepackage{CJK}
\usepackage{epsfig}
\usepackage{amsthm}
\usepackage{amsmath}
\usepackage{amsbsy}
\usepackage{amsfonts}
\usepackage{graphicx}
\usepackage{algorithm}
\usepackage{algorithmic}
\usepackage{rotating}
\usepackage{enumerate}
\usepackage{float}
\newtheorem{theorem}{Theorem}[section]

\theoremstyle{definition}

\theoremstyle{remark}

\numberwithin{equation}{section}
\begin{document}

\title{A Residual Based Sparse Approximate Inverse Preconditioning
Procedure for Large Sparse Linear Systems\thanks{Supported in part by
National Science Foundation of China (No. 11371219).}}
\author{Zhongxiao Jia\thanks{Department of Mathematical
Sciences, Tsinghua University, Beijing, 100084,
People's Republic of China, jiazx@tsinghua.edu.cn.}
\and Wenjie Kang\thanks{Department of Mathematical Sciences,
Tsinghua University, Beijing, 100084,
People's Republic of China, kangwj11@mails.tsinghua.edu.cn.}}

\date{}
\maketitle
\begin{abstract}
The SPAI algorithm, a sparse approximate inverse preconditioning
technique for large sparse linear systems, proposed by Grote and Huckle
[SIAM J. Sci. Comput., 18 (1997), pp.~838--853.],
is based on the F-norm minimization and computes a sparse
approximate inverse $M$ of a large sparse matrix $A$ adaptively.
However, SPAI may be costly to seek the most profitable indices at each loop
and $M$ may be ineffective for preconditioning. In this paper,
we propose a residual based sparse approximate inverse preconditioning
procedure (RSAI), which, unlike SPAI, is based on only the {\em dominant} rather
than all information on the current residual and augments sparsity patterns
adaptively during the loops. RSAI is less costly to
seek indices and is more effective to capture a good approximate
sparsity pattern of $A^{-1}$ than SPAI. To control the sparsity of
$M$ and reduce computational cost, we develop a practical RSAI($tol$)
algorithm that drops small nonzero entries adaptively during the process.
Numerical experiments are reported to demonstrate that RSAI($tol$)
is at least competitive with SPAI and can be considerably more efficient
and effective than SPAI. They also indicate that RSAI($tol$) is comparable to
the PSAI($tol$) algorithm proposed by one of the authors in
2009.
\smallskip

\smallskip


\end{abstract}

\section{Introduction}\label{sec:1}
Consider the iterative solution of large sparse linear system
\begin{equation}\label{eq1}
  Ax=b,
\end{equation}
where $A$ is an $n\times n$ real nonsingular and nonsymmetric
matrix, and $b$ is a given $n$-dimensional vector.
This kind of problem is a core problem in scientific and engineering computing.
Krylov iterative solvers, such
as the generalized minimal residual method (GMRES) % \cite{saad86},
and the biconjugate gradient stabilized method (BiCGStab) \cite{barrett94,saad03},
have been commonly used in nowdays for solving \eqref{eq1}. % \cite{vorst92}.
However, when $A$ has bad spectral property or is ill conditioned,
the convergence of Krylov solvers are generally extremely
slow \cite{saad03}. In order to accelerate the convergence of Krylov solvers,
one must utilize
preconditioning techniques to improve the conditioning of \eqref{eq1},
so that Krylov solvers applied to resulting preconditioned systems converge fast.
Sparse approximate inverse (SAI) preconditioning has been one major class of
general-purpose preconditioning \cite{benzi02,Ferronato12,saad03},
and it aims at computing a sparse approximate inverse
$M\approx A^{-1}$ or factorized $M=M_1M_2\approx A^{-1}$ directly. With such
an $M$ available, the right and left preconditioned systems are
\begin{equation}\label{eq2}
AMy=b,x=My \quad\text{and}\quad MAx=Mb,
\end{equation}
respectively, and the factorized preconditioned system is
\begin{equation}\label{eq3}
M_1AM_2y=M_1b,\ x=M_2y.
\end{equation}
We then use a Krylov solver to solve \eqref{eq2}
or \eqref{eq3}, depending on the way that $M$ is applied. Since the
coefficient matrices in the above preconditioned systems are roughly the
identity matrix $I$, Krylov solvers are expected to converge quickly.

The success of SAI preconditioning is based on the underlying
hypothesis that the majority of entries of $A^{-1}$ are small,
which means that $A$, indeed, has good sparse approximate inverses.
A good preconditioner $M$ should be as sparse as possible,
and it should be constructed efficiently and applied within Krylov solvers
cheaply. There are two kinds of approaches to computing $M$. One of them gets
a factorized $M=M_1M_2$ and applies $M_1$ and $M_2$ to \eqref{eq3}.
Efficient algorithms of this kind are approximate inverse (AINV) type
algorithms, which are
derived from the incomplete biconjugation procedure \cite{benzi96,benzi98}.
Their stabilized and block variations are developed in \cite{benzi01}.
An alternative is the balanced incomplete factorization (BIF) algorithm,
which computes an incomplete LU (ILU) factorization
and its inverse simultaneously \cite{tuma08,tuma10}.
The other kind of approach is based on F-norm minimization, which is
inherently parallelizable and constructs $M$ by
minimizing $\|AM-I\|_{F}$ with certain sparsity constraints
on $M$, where $\|\cdot\|_{F}$ denotes the Frobenius norm of a matrix.
We will introduce more on this approach in the next paragraph.
Kolotilina and Yeremin \cite{kolotilina93} have proposed a factorized
sparse approximate inverse (FSAI) preconditioning procedure, which is a mixture
of the above two kinds. FSAI has been generalized to
block form, called BFSAI, in \cite{Janna10,Janna11}. An adaptive
algorithm that generates the pattern of the BFSAI
preconditioner $M$ can be found in \cite{Janna14,Janna11,Janna13}.
For a comprehensive survey and comparison of SAI preconditioning procedures,
we refer the reader to \cite{benzi02,benzi99}.

We focus on F-norm minimization based SAI preconditioning and revisit
the SPAI algorithm \cite{Grote97} in this paper.
A key of this kind of preconditioning is
the determination of an effective approximate sparsity pattern of $A^{-1}$.
There are two approaches to doing this, one of which is static and the other
is adaptive. A static SAI preconditioning procedure first prescribes a
sparsity pattern of $M$ and then computes $M$ by solving $n$ least
squares (LS) problems independently \cite{benson82,benson84}.
The main difficulty of this approach is how to choose an effective approximate
sparsity pattern of $A^{-1}$. A lot of research has
been done on this issue, and some priori envelope patterns for effective
approximate sparsity patterns of $A^{-1}$  have been established;
see, e.g., \cite{chow00,gilbert94,Huckle99}.
For a general irreducible sparse $A$, however, these envelope patterns are
often quite dense, so that it is expensive to use them as patterns
of $M$ directly. To this end,
several researchers have proposed and developed adaptive procedures,
which start with a simple initial sparsity pattern and successively augment
or adjust it until either the resulting $M$ satisfies a prescribed accuracy,
or the maximum loops are performed, or
the maximum number of nonzero entries in $M$ is reached.
Such idea was first advocated in \cite{cosgrove92}. Grote and Huckle
\cite{Grote97} have proposed the SPAI algorithm aimed at augmenting the
sparsity pattern of $M$ adaptively by adding the small number of most profitable
indices at each loop. It has been generalized
to block form, called BSPAI, in \cite{barnard99}.
Gould and Scott \cite{gould98} have given a number of enhancements which
may improve the performance of the SPAI algorithm.
Chow and Saad \cite{chow98} have put forward a minimal
residual based (MR) algorithm that uses the sparse-sparse iteration with
dropping strategies. SPAI is more robust than the MR algorithm \cite{benzi02}.
Motivated by the Cayley--Hamilton theorem,
Jia and Zhu \cite{Jia09} have proposed an adaptive Power SAI
(PSAI) preconditioning procedure and developed a practical
PSAI($tol$) algorithm with $tol$ dropping tolerance, which has been
shown to be at least competitive with SPAI and can outperform
SPAI considerably for some difficult problems.
Jia and Zhang \cite{Jia13a} have recently established a mathematical
theory on dropping tolerances $tol$ for PSAI($tol$) and all
the static F-norm minimization based SAI preconditioning procedures.
Based on the theory, they have designed robust adaptive dropping criteria.
With the criteria applied, PSAI($tol$) and the mentioned static
SAI preconditioning procedures can make $M$ as sparse as possible and as equally
effective as the possibly much denser one generated by the basic PSAI or the
static SAI procedures without dropping small entries.

Remarkably, the unique fundamental mathematical distinction of all the adaptive
F-norm minimization based SAI preconditioning procedures consists in the way
that the sparsity pattern of $M$ is augmented or adjusted.
Practically, both static and adaptive SAI procedures must control the sparsity
of $M$. We mention that there is a
prefiltration, which shrinks the pattern of $A$ by dropping small entries of $A$
before implementing a SAI preconditioning
algorithm \cite{Kaporin92,Vavasis92,chen98,tang99}.

Jia and Zhang \cite{Jia13b} have made an analysis on SPAI and
PSAI($tol$) and shown why PSAI($tol$) can be more effective than SPAI for
preconditioning \eqref{eq1}, accompanied by detailed numerical comparisons of
the two algorithms on a lot of regular and irregular sparse problems arising
from applications. Here the meaning of `regular sparse' is that all columns
of $A$ are sparse, and that of `irregular sparse' is that $A$ has at least
one relatively dense column, whose number of nonzero entries is
substantially more than the average number of nonzero entries per column of $A$.
Empirically and numerically, a column is declared irregular sparse
if it has at least $10p$
nonzero entries, where $p$ is the average number of nonzero entries per column
of $A$ \cite{Jia13b}. For $A$ irregular sparse,
Jia and Zhang \cite{Jia13b} have shown that SPAI must be costly to seek and add
indices at each loop and, moreover, the resulting $M$ may be ineffective
for preconditioning, while PSAI($tol$), though also very costly,
can produce an effective
preconditioner $M$. To this end, they have proposed an approach that first
transforms the irregular sparse \eqref{eq1} into certain new regular ones
and then uses SPAI and PSAI($tol$) to construct $M$ for the regular
problems. Such approach greatly improves the computational efficiency
of SPAI and PSAI($tol$) as well as the preconditioning effectiveness of SPAI
applied to irregular sparse \eqref{eq1} directly.

In this paper, suppose that the current $M$ is not good enough, and we
need to augment or adjust the sparsity pattern of $M$ in order to
get a better $M$. We will propose a new adaptive
F-norm minimization based SAI preconditioning procedure,
called the residual based SAI (RSAI) algorithm, and develop
a practical RSAI algorithm with dropping criteria exploited.
The RSAI algorithm may greatly improve the computational efficiency
and the preconditioning effectiveness of SPAI, especially when $A$ is
irregular sparse. The basic idea of RSAI algorithm is as follows: Differently
from SPAI, for each column $m_k$, $k=1,2,\ldots,n$, of the current $M$,
by only selecting a few dominant indices that correspond
to the largest entries in the residual of $m_k$,
we augment the sparsity pattern of $m_k$ using a new approach that
avoids some possibly expensive computation and logical comparisons in SPAI
when determining the most profitable indices based on the whole residual of
$m_k$. As it will turn out, the RSAI procedure is not only (considerably)
less costly to seek and add indices but also more effective to
capture a good approximate sparsity pattern of $A^{-1}$ than SPAI,
which is especially true for an irregular sparse $A$.
We derive a quantitative estimate for the number of nonzero entries
in $M$, demonstrating how it depends on the sparsity pattern of $A$,
the number of indices exploited that correspond to the largest
entries of the residual at each loop, and the number $l_{\max}$ of loops.
To control the sparsity of $M$ and improve computational efficiency,
we develop a practical RSAI($tol$) algorithm by making use of the
adaptive dropping criterion established in \cite{Jia13a},
which guarantees that two $M$ obtained by RSAI($tol$) and the basic RSAI
without dropping small entries have comparable preconditioning effects.
We show that the positions of large entries in $M$ are automatically adjusted
in a global sense during the loops. It is known \cite{Jia13a,Jia09}
that SPAI retains the already occupied
positions of nonzero entries in $M$ in subsequent loops and adds new positions
of nonzero entries in $M$ at each loop. As a result,
the RSAI($tol$) algorithm captures an approximate sparsity pattern of $A^{-1}$
in a globally optimal sense, while the SPAI algorithm achieves this goal only
in a locally optimal sense. This difference may make RSAI($tol$) advantageous
over SPAI. Numerical experiments will confirm that
RSAI($tol$) is at least competitive with and can be
substantially more efficient and effective than SPAI, and they will also
illustrate that RSAI($tol$) is as comparably effective as PSAI($tol$).

The paper is organized as follows. In Section \ref{sec:2},
we review the F-norm minimization based SAI preconditioning and
the SPAI procedure, and introduce the notation to be used.
In Section \ref{sec:3}, we propose the basic RSAI procedure.
In Section \ref{sec:4}, we make a theoretical analysis and give
some practical considerations, based on which we develop a practical
RSAI($tol$) algorithm with dynamic dropping strategy in \cite{Jia13a}
exploited. Finally, we make numerical experiments on a
number of real-world problems to confirm our assertions
on RSAI($tol)$, SPAI and PSAI($tol$) in Section \ref{sec:5}.

\section{The F-norm minimization SAI preconditioning and
the SPAI procedure} \label{sec:2}

A F-norm minimization based SAI preconditioning procedure
solves the problem
\begin{equation}\label{minfnorm}
\min_{M\in \mathcal{M}}\|AM-I\|_{F},
\end{equation}
where $\mathcal{M}$ is the set of matrices with a given sparsity
pattern $\mathcal{J}$. Denote by $\mathcal{M}_{k}$ the set of
$n$-dimensional vectors whose sparsity pattern is
$\mathcal{J}_{k}=\{i\mid(i,k)\in \mathcal{J}\}$, and let
$M=(m_1,m_2,\ldots,m_n)$. Then
\eqref{minfnorm} can be separated into $n$ independent
constrained LS problems
\begin{equation}\label{min2norm}
  \min_{m_{k}\in \mathcal{M}_{k}}\|Am_{k}-e_{k}\|,\ k=1,2,\ldots,n,
\end{equation}
where $\|\cdot\|$ denotes the 2-norm of a matrix or vector,
and $e_{k}$ is the $k$-th column of the identity
matrix $I$ of order $n$. For each $k$, let $\mathcal{I}_{k}$ be the set of
indices of nonzero rows of $A(\cdot,\mathcal{J}_{k})$.
Define $A_k=A(\mathcal{I}_{k},\mathcal{J}_{k})$, the
reduced size vector $\tilde{m}_k=m_{k}(\mathcal{J}_{k})$, and
$\tilde{e}_{k}=e_{k}(\mathcal{I}_{k})$. Then
\eqref{min2norm} amounts to solving the smaller unconstrained LS problems
\begin{equation}\label{reducemin2norm}
\min_{\tilde{m}_{k}}\|A_k\tilde{m}_{k}-\tilde{e}_{k}\|, k=1,2,\ldots,n,
\end{equation}
which can be solved by QR decompositions in parallel.

If $M$ is not good enough, an adaptive SAI preconditioning procedure,
such as SPAI \cite{Grote97} and PSAI \cite{Jia09},
improves it by augmenting or adjusting the sparsity
pattern $\mathcal{J}_{k}$ dynamically and updating
$\tilde{m}_{k}$ for $k=1,2,\ldots,n$ efficiently. We describe SPAI below.

Denote by $\mathcal{J}_{k}^{(l)}$ the sparsity pattern of $m_{k}$
after $l$ loops starting with an initial
pattern $\mathcal{J}_{k}^{(0)}$, and by $\mathcal{I}_{k}^{(l)}$
the set of indices of nonzero rows of $A(\cdot,\mathcal{J}_{k}^{(l)})$.
Let $A_k=A(\mathcal{I}_{k}^{(l)},\mathcal{J}_{k}^{(l)})$, $\tilde{e}_k
=e_{k}(\mathcal{I}_{k}^{(l)})$ and $\tilde{m}_{k}$ the solution of
\eqref{reducemin2norm}. Denote the residual of \eqref{min2norm} by
\begin{equation}\label{residual}
r_k=A m_k-e_k,
\end{equation}
whose norm is exactly equal to the residual norm $\|\tilde{r}_k\|$
of \eqref{reducemin2norm} defined by $\tilde{r}_{k}=A_k\tilde{m}_{k}
-\tilde{e}_{k}$.
If $r_{k}\neq0$, define $\mathcal{L}_{k}$ to be the set of indices $i$
for which $r_{k}(i)\neq 0$ and $\mathcal{N}_{k}$ the set
of indices of nonzero columns of $A(\mathcal{L}_{k},\cdot)$. Then
\begin{equation}\label{argument}
  \hat{\mathcal{J}}_{k}=\mathcal{N}_{k}\setminus\mathcal{J}_{k}^{(l)}
\end{equation}
constitutes the new candidates for augmenting $\mathcal{J}_{k}^{(l)}$
in the next loop of SPAI. For each $j\in \hat{\mathcal{J}}_{k}$, SPAI solves the
one-dimensional problem
\begin{equation}\label{onedimensional}
  \min_{\mu_{j}}\|r_{k}+\mu_{j}Ae_{j}\|,
\end{equation}
and the 2-norm $\rho_{j}$ of the new residual $r_{k}+\mu_{j}Ae_{j}$ satisfies
\begin{equation}\label{2normofresidual}
  \rho_{j}^{2}=\|r_{k}\|^{2}-\frac{(r_{k}^{T}Ae_{j})^{2}}
  {\|Ae_{j}\|^{2}}.
\end{equation}
SPAI takes $l_a$, $1\sim 5$, indices from $\hat{\mathcal{J}}_{k}$
corresponding to the smallest $\rho_{j}$, called the most profitable indices,
and adds them to $\mathcal{J}_{k}^{(l)}$
to obtain $\mathcal{J}_{k}^{(l+1)}$.
Let $\hat{\mathcal{I}}_{k}$ be the set of indices
of new nonzero rows corresponding to the most profitable indices added.
SPAI gets the new row indices $\mathcal{I}_{k}^{(l+1)}=\mathcal{I}_{k}^{(l)}
\bigcup \hat{\mathcal{I}}_{k}$, and by it and $\mathcal{J}_{k}^{(l+1)}$
we form an updated LS problem \eqref{reducemin2norm}, which is cheaply
solved by updating the previous $\tilde{m}_{k}$. Proceed in such way until
$\|r_k\|=\|Am_{k}-e_{k}\|\leq\varepsilon$
or $l\geq l_{\max}$, where
$\varepsilon$ is a prescribed tolerance, usually $0.1\sim0.4$.

Each loop of SPAI consists of two main steps, which
include the selection of the most profitable indices and
the solution of the resulting new LS problem, respectively.
For the selection of the most profitable indices,
one first determines $\mathcal{L}_{k}$
through $r_{k}$ and $\hat{\mathcal{J}}_{k}$ through $\mathcal{L}_{k}$ and
$\mathcal{J}_{k}$, computes the $\rho_j$, then orders them, and finally
selects the most profitable indices. Clearly, whenever the
cardinality of $\hat{\mathcal{J}}_{k}$ is big, this step is
time consuming. It has been shown \cite{Jia13b}
that the cardinality of $\hat{\mathcal{J}}_{k}$ is
always big for $l=1,2,\ldots,l_{\max}$ when the $k$-th column of $A$ is
relatively dense and the initial pattern $\mathcal{J}^{(0)}$ is that of
$I$, causing that SPAI is very costly to select the most profitable indices.
In addition, we can easily see that SPAI is also costly to seek the most 
profitable indices for $A$ row irregular sparse. This is the case 
that an index in $\mathcal{L}_k$ corresponds to a relatively dense row 
of $A$. For detailed numerical evidence on this, we refer
to \cite{Jia13b,jia15}, where it has been clearly shown that SPAI 
is too costly and impractical since it spends too
much time seeking the most profitable indices when $A$ is irregular sparse.

\section{The RSAI preconditioning procedure}\label{sec:3}

Our motivation for proposing a new SAI preconditioning procedure
is that SPAI may be costly to select the most profitable indices
and may be ineffective for preconditioning, which is definitely
the case when $A$ has some relatively dense columns.
Our new approach to augmenting $\mathcal{J}_{k}^{(l)}$ and
$\mathcal{I}_{k}^{(l)}$
is based on the partially {\em dominant} other than all information
on the current residual $r_k$, and picks up new indices
at each loop directly and cheaply that avoids possibly expensive
computation and logical comparisons needed by SPAI. Importantly,
the new SAI preconditioning procedure is more effective to capture a good
approximate sparsity pattern of $A^{-1}$. Since our
procedure critically depends on the {\em sizes} of entries of $r_k$, we call
the resulting procedure the Residual based Sparse Approximate
Inverse (RSAI) preconditioning procedure. In what follows we present
a basic RSAI procedure.

Suppose that $M$ is the one generated by RSAI after $l$
loops starting with the initial sparsity pattern $\mathcal{J}^{(0)}$.
Consider the residual $\|r_k\|$
defined by \eqref{residual} for $k=1,2,\ldots,n$.
If $\|r_k\|\leq \varepsilon$ for a prescribed tolerance $\varepsilon$,
then $m_{k}$ satisfies the required accuracy, and we do not
improve $m_k$ further. If $\|r_{k}\|>\varepsilon$, we must augment or
adjust $\mathcal{J}_{k}^{(l)}$ to update $m_k$ so as
to reduce $\|r_{k}\|$.

Denote by $r_{k}(i)$ the $i$-th entry of $r_k$,
and by $\mathcal{L}_{k}$ the set of indices $i$ for which
$r_{k}(i)\neq0$.
{\em Heuristically}, the indices corresponding to the largest entries
$r_{k}(i)$ of $r_{k}$ in magnitude are the most important and
dominate $\mathcal{L}_{k}$ in the sense that these large entries contribute
most to the size of $\|r_k\|$. Therefore, in order to reduce $\|r_k\|$
both substantially and cheaply, these most important or dominant
indices should have priority, that is,
we should take precedence to reduce the large entries $r_{k}(i)$ of $r_k$
by augmenting or adjusting the pattern of $m_k$ based on the dominant
indices described above.
Therefore, unlike SPAI, as a starting point, rather than using the whole
$\mathcal{L}_k$, RSAI will exploit only the dominant subset of it
and augment or adjust the pattern
of $m_k$ in a completely new manner, as will be detailed below.

At loop $l$, denote by $\hat{\mathcal{R}}_{k}^{(l)}$ the set of the dominant
indices $i$ corresponding to the largest $|r_{k}(i)|$.
Let $\hat{\mathcal{J}}_{k}$ be
the set of all new column indices of $A$ that correspond to
$\hat{\mathcal{R}}_{k}^{(l)}$ of row indices
but do not appear in $\mathcal{J}_{k}^{(l)}$.
We then add $\hat{\mathcal{J}}_{k}$ to
$\mathcal{J}_{k}^{(l)}$  to obtain a new sparsity pattern
$\mathcal{J}_{k}^{(l+1)}$.
Denote by $\hat{\mathcal{I}}_{k}$
the set of indices of new nonzero rows corresponding to the added
column indices $\hat{\mathcal{J}}_{k}$. Then we update
$\mathcal{I}_{k}^{(l+1)}=\mathcal{I}_{k}^{(l)}\bigcup\hat{\mathcal{I}}_{k}$
and form the new LS problem
\begin{equation}\label{upmin2norm}
    \min\|A(\mathcal{I}_{k}^{(l+1)},\mathcal{J}_{k}^{(l+1)})m_{k}
    (\mathcal{J}_{k}^{(l+1)})-e_{k}(\mathcal{I}_{k}^{(l+1)})\|,
\end{equation}
whose solution can be updated from $m_{k}(\mathcal{J}_{k}^{(l)})$
in the way described in \cite{Grote97,Jia09}, and the
updated $m_k$ is a better approximation
to the $k$-th column of $A^{-1}$. We repeat this process until
$\|r_k\|\leq\varepsilon$ or $l$ reaches $l_{\max}$.
The above RSAI procedure effectively suppresses the effects
of the large entries $|r_{k}(i)|$ and reduces $\|r_k\|$.

Now we give more insight into $\hat{\mathcal{R}}_{k}^{(l)}$. When choosing it
from $\mathcal{L}_{k}$ in the above way, we may encounter
$\hat{\mathcal{R}}_{k}^{(l)}=\hat{\mathcal{R}}_{k}^{(l-1)}$. If so,
we cannot augment the sparsity pattern of $m_k$. In this case, we set
\begin{equation}\label{insertR}
  \mathcal{R}_{k}^{(l)}=\bigcup_{i=0}^{l-1}
  \hat{\mathcal{R}}_{k}^{(i)},
\end{equation}
and choose $\hat{\mathcal{R}}_{k}^{(l)}$ from the set whose elements are in
$\mathcal{L}_{k}$ but not in $\mathcal{R}_{k}^{(l)}$. Obviously, the resulting
$\hat{\mathcal{R}}_{k}^{(l)}$ is always non-empty unless $m_k$ is exactly
the $k$-th column of $A^{-1}$.
On the other hand, if $\hat{\mathcal{J}}_{k}$ happens to be empty,
then $\mathcal{J}_{k}^{(l)}=\mathcal{J}_{k}^{(l+1)}$ and we just skip
to loop $l+2$, and so forth. Since $r_{k}\neq0$, there must exist
a $\check{l}\geq l+1$ such that $\hat{\mathcal{J}}_{k}$ is not empty.
Otherwise, $\mathcal{J}_{k}^{(l)}$ is the set of indices of all nonzero
columns of $A(\mathcal{L}_{k},\cdot)$, and
\begin{equation}\label{keepfillin}
r_{k}(\mathcal{L}_{k})^{T}A(\mathcal{L}_{k},\mathcal{J}_{k}^{(l)})
 = r_{k}(\mathcal{L}_{k})^{T}A(\mathcal{L}_{k},\cdot)=0,
\end{equation}
which means that $r_{k}(\mathcal{L}_{k})=0$, i.e., $r_k=0$,
and $m_k$ is exactly the $k$-th column of $A^{-1}$
since $A(\mathcal{L}_{k},\cdot)$ has row full rank.

The above RSAI procedure can be described as Algorithm 1, named as the
basic RSAI algorithm.
\begin{algorithm}[htb]
\caption{The basic RSAI Algorithm}\label{RSAI}
\vspace{.1cm}
 For $k=1,2,\ldots,n$ compute $m_{k}$:
\begin{algorithmic}[1]\label{nalg}
\vspace{.2cm}
\STATE Choose an initial pattern $\mathcal{J}_{k}^{(0)}$ of $m_{k}$, and give a
user-prescribed tolerance $\varepsilon$ and the maximum number $l_{\max}$
of loops. Set $l=0$ and $\mathcal{R}_{k}^{(0)}=\emptyset$.

\STATE Determine $\mathcal{I}_{k}^{(0)}$, the set of indices of nonzero
rows of $A(\cdot,\mathcal{J}_{k}^{(0)})$,
solve \eqref{reducemin2norm} for $\tilde{m}_k$
by the QR decomposition of $A_k=A(\mathcal{I}_{k}^{(0)},
\mathcal{J}_{k}^{(0)})$, recover $m_{k}$ from $\tilde{m}_k$,
and compute $r_{k}=Am_{k}-e_{k}$.

\WHILE{$\|r_k\|>\varepsilon$ and $l < l_{\max}$}

\STATE Set $\mathcal{L}_{k}$ to be the set of indices $i$ for which
$r_{k}(i)\neq 0$, sort $|r_{k}(i)|$ in decreasing order,
and let $\hat{\mathcal{R}}_{k}^{(l)}$ be
the set of dominant indices $i$ that correspond to a few largest
$|r_k(i)|$ appearing in $\mathcal{L}_{k}$ but not
in $\mathcal{R}_{k}^{(l)}$.
Augment
$\mathcal{R}_{k}^{(l+1)}=\mathcal{R}_{k}^{(l)}\bigcup\hat{\mathcal{R}}_{k}^{(l)}$.

\STATE Set $\hat{\mathcal{J}}_{k}$ equal to the set of all new column
indices of $A(\hat{\mathcal{R}}_{k},\cdot)$ but not in $\mathcal{J}_{k}^{(l)}$.
Let $\mathcal{\hat{I}}_k$ be the set of row indices corresponding to
$\mathcal{\hat{J}}_k$ that do not appear in $\mathcal{I}_k^{(l)}$, and update
$\mathcal{I}_k^{(l+1)}=\mathcal{I}_k^{(l)}\bigcup \mathcal{\hat{I}}_k$.

\STATE Set $l=l+1$; if $\hat{\mathcal{J}}_{k}=\emptyset$, then go to step 3.

\STATE For each $j\in\mathcal{\hat{J}}_k$, update $m_{k}$ and $\|r_k\|$ using the
approach in \cite{Grote97,Jia09}, respectively.

%\STATE If $\|r_k\|\leq\varepsilon$, break;
%else go to step 3.

\ENDWHILE
\end{algorithmic}
\end{algorithm}

Two key differences between SPAI and RSAI are clear now. Firstly,
for RSAI we order the nonzero entries in $r_k$, pick up a few
dominant indices $\mathcal{\hat{R}}_k^{(l)}$
with the largest entries of $r_k$
in magnitude that do not appear in $\mathcal{R}_k^{(l)}$.
Using $\mathcal{\hat{R}}_k^{(l)}$, we
determine the new indices $\mathcal{\hat{J}}_k$ and add them to
$\mathcal{J}_k^{(l)}$
to get the new column indices $\mathcal{J}_k^{(l+1)}$, by which we identify
the new row indices $\hat{\mathcal{I}}_k$ to be added and form the
set $\mathcal{I}_k^{(l+1)}$ of row indices. Secondly, for RSAI we do not
perform possibly expensive steps
\eqref{onedimensional} and \eqref{2normofresidual} and the ordering and
sorting followed. Since $\mathcal{\hat{R}}_k^{(l)}$
is a subset of $\mathcal{L}_k$ and we assume that it has only
a few elements, its cardinality can be much smaller
than that of $\mathcal{L}_k$, which is typically true for $A$ irregular
sparse. As a result, the determination of $\mathcal{I}_{k}^{(l+1)}$
and $\mathcal{J}_{k}^{(l+1)}$ is less costly, and it
can be substantially less time consuming than SPAI,
especially when $A$ is irregular sparse.

Similar to SPAI, we need to provide an initial sparsity pattern
of $M$ for the RSAI algorithm, which is usually chosen to be
that of $I$ when $A$ has nonzero diagonals. We also need to provide
a stopping criterion $\varepsilon$, the number of the dominant
indices and the maximum number $l_{\max}$ of loops. For them, we
take $\varepsilon=0.1\sim 0.4$, similarly to that used in F-norm based
SAI preconditioning procedures including SPAI and PSAI. We suggest taking
the cardinality $c$ of $\mathcal{\hat{R}}_k^{(l)}$ to be 3 or so
at each loop. As for outer loops $l_{\max}$, we take it to be small, say 10.

In the next section, we make some theoretical
analysis and develop a more practical RSAI algorithm with some dropping
strategy used.

\section{Theoretical analysis and a practical RSAI algorithm}\label{sec:4}

We cannot guarantee that $M$ obtained by the basic RSAI
algorithm is nonsingular without additional
requirements. Grote and Huckle \cite{Grote97} present several results,
showing how the non-singularity of $M$ is related to $\varepsilon$ and
how the eigenvalues and the singular values of the preconditioned matrix
$AM$ distribute. Their results are general and apply to
$M$ obtained by any F-norm minimization based SAI preconditioning
procedure.

Huckle \cite{Huckle99} shows that the patterns of $(A^{T}A)^{\mu-1}A^{T}$
for small $\mu$ are effective upper bounds for the sparsity pattern of $M$
by the SPAI algorithm. Note that both RSAI and SPAI augment the sparsity
pattern of $M$ based on the indices of nonzero entries of
residuals $r_k$, $k=1,2,\ldots,n$. Consequently,
in the same way as \cite{Huckle99}, we can justify that the patterns of
$(A^{T}A)^{\mu-1}A^{T}$ for small $\mu$ are also upper bounds for the
sparsity pattern of $M$ obtained by the basic RSAI algorithm.

Let us get insight into the pattern of $M$ by the basic RSAI algorithm.
Obviously, the dominant indices at each loop and the maximum
number $l_{\max}$ of loops directly affect the sparsity of
$M$. Now we present a quantitative upper
bound for the number of nonzero entries in $M$ and show how the sparsity
of $M$ depends on that of $A$, the number $c$ of the dominant indices
at each loop and $l_{\max}$.

\begin{theorem}\label{theorem1}
Assume that RSAI runs $l_{\max}$ loops to generate $M=(m_1,m_2,\ldots,m_n)$.
Denote by $g_{k}$ the number of nonzero entries in $A(k,:)$, $1\leq k\leq n$, by
$g=\max_{1\leq k\leq n}{g_{k}}$, by $c$ the number of the
dominant indices at each loop, and by $nnz(m_k)$
the number of nonzero entries in $m_{k}$.
Then for $\mathcal{J}_{k}^{(0)}=\{k\}$ we have
\begin{equation}\label{estimatenonzeromk}
  nnz(m_k)\leq\min\{g c l_{\max}+1,n\}, k=1,2,\ldots,n
\end{equation}
and
\begin{equation}\label{estimatenonzerom}
 nnz(M)\leq\min\{(g c l_{\max}+1)n,n^{2}\}.
\end{equation}
\end{theorem}

\begin{proof}
By assumption, it is known that the number of nonzero entries added
to $m_{k}$ at each loop does not exceed $g c$. Therefore, we have
\begin{equation*}
 nnz(m_k)\leq g c l_{\max}+1,
\end{equation*}
which, together with the trivial bound
$nnz(m_k)\leq n$, establishes \eqref{estimatenonzeromk}.
Note that the number of nonzero entries of $M$ is $\sum_{k=1}^n nzz(m_k)$.
\eqref{estimatenonzerom} is direct from \eqref{estimatenonzeromk}.
\end{proof}

Theorem \ref{theorem1} shows that if all the rows of $A$ are sparse,
i.e., $g$ is small, then $M$ must be sparse for $l_{\max}$ and $m$ small.
On the other hand, if $A$ has one relatively dense row, some columns of
$M$ may become dense quickly with increasing $l$. This is the case
once an index in $\mathcal{\hat{R}}_k^{(l)}$ at some loops
corresponds to a relatively dense row of $A$. In this case, the basic
RSAI is inefficient since a relatively large LS problem will emerge
in the next loop, causing that solving it is expensive.
This is a shortcoming of the basic RSAI algorithm for $A$ row irregular
sparse.

Under the underlying hypothesis that the majority of the
entries of $A^{-1}$ are small, we know that
whenever $M$ becomes relatively dense in the
course of construction, the majority of its entries must be small and make
little contribution to $A^{-1}$. Therefore, in order to control the
sparsity of $M$ and improve the efficiency of the basic RSAI algorithm,
we should drop those small entries promptly during the process for
practical purposes. This asks us to introduce some reasonable dropping
strategies into the basic RSAI algorithm so as to develop practical RSAI
algorithms for constructing an effective and sparse $M$.

We should be aware that SPAI implicitly uses a dropping strategy to
ensure that all the columns of $M$ constructed by it are sparse.
Precisely, suppose that SPAI is run $l_{\max}$ loops starting with
the pattern of $I$ and a few, say, $l_a$ most profitable indices are
added to the pattern of $m_k$ at each loop. Then the number of nonzero
entries of the final $m_k$ does not exceed $1+ l_al_{\max}$,
which is fixed and small as $l_a$ and $l_{\max}$ are
both fixed and small. Such $M$ may not be robust since the number of large
entries in the $k$-th column of $A^{-1}$ is unknown in practice.
Moreover, for a general sparse $A$, the numbers of large entries in
the columns of $A^{-1}$ may have great differences, that is, good sparse
approximations of $A^{-1}$ may well be irregular sparse. This is particularly
true for $A$ irregular sparse and even for $A$ regular sparse \cite{Jia13b}.
Consequently, SPAI may generate a poor preconditioner $M$ since
some columns of it are too sparse for given small $l_a$ and $l_{\max}$.

The above analysis suggests that we should not fix the number of large
entries of each column of $M$ for RSAI in advance.
In the spirit of PSAI($tol$) \cite{Jia09},
a more robust and general-purpose dropping strategy is to retain all
the large entries of $m_k$ produced and drop those small ones below
a prescribed tolerance $tol$ during the loops.
This kind of dropping strategy should better capture a good approximate sparsity
pattern of $A^{-1}$ and produce an effective $M$ more possibly.

Dropping tolerances $tol$ used in SAI preconditioning procedures had been
empirically chosen as some small quantities, say $10^{-3}$.
%Some research has
%been done on this issue, see, e.g., \cite{bollhofer02,chow98,Janna11}.
Recently, Jia and Zhang \cite{Jia13a} have shown that such
dropping criteria are not robust and may lead
to a sparser but ineffective preconditioner $M$ or
a possibly quite dense $M$, which, though effective for
preconditioning, is very costly to construct and apply. Jia and
Zhang \cite{Jia13a} have established a mathematical theory on robust dropping
tolerances for PSAI($tol$) and all the static F-norm minimization based SAI
preconditioning procedures. Based on the theory, they have designed robust and
adaptive selection criteria for dropping tolerances, which
adapt to the RSAI algorithm directly: an entry $m_{jk}$ is dropped whenever
\begin{equation}\label{dropping}
  |m_{jk}|\leq tol_{k}=\frac{\varepsilon}{nnz(m_{k})\|A\|_{1}},j=1,2,\ldots,n
\end{equation}
during the loops, where $m_{jk}$ is the $j$-th entry of $m_{k}$ and $nnz(m_{k})$
is the number of nonzero entries in $m_{k}$ at the current loop,
and $\|\cdot\|_1$ is the 1-norm of a matrix. In terms of the theory
in \cite{Jia13a}, the RSAI($tol$) equipped with the above dropping criterion
will generate a SAI preconditioner $M$ that has comparable preconditioning
quality to the possibly much denser one generated by the basic RSAI algorithm
without dropping small nonzero entries; see \cite[Theorem 3.5]{Jia13a}.

Introducing \eqref{dropping} into Algorithm \ref{RSAI},
we have developed a practical algorithm, called the RSAI($tol$)
algorithm. As a key comparison of RSAI($tol$) and SPAI,
we notice that, for RSAI($tol$), the positions of large entries of
$m_{k}$ are adjusted dynamically during the loops, while
the already occupied positions of
nonzero entries in $m_{k}$ by SPAI retain unchanged in
subsequent loops and we simply add a few new indices to
the pattern of $m_{k}$ at each loop. Remarkably,
some entries of $m_k$ are well likely to change from large to small
during the loops, so that the final $M$ may have some small entries
that contribute little to $A^{-1}$. In other words, RSAI($tol$)
seeks the positions of large entries of $A^{-1}$ in a globally optimal
sense, while the SPAI algorithm achieves this goal in a locally optimal sense.
Consequently, RSAI($tol$) captures the
sparsity pattern of $A^{-1}$ more effectively than SPAI.

\section{Numerical experiments}\label{sec:5}

In this section we test a number of real-world problems coming from
applications, which are described in
Table \ref{table-mtr}\footnote{All of these matrices are either from the
Matrix Market of the National Institute of Standards and Technology at
http://math.nist.gov/MatrixMarket/ or from the University of Florida Sparse
Matrix Collection at http://www.cise.ufl.edu/research/sparse/matrices/}.
We also list some useful information about the test matrices in Table
\ref{table-mtr}. For each matrix, we give the number $s$ of irregular columns
as defined in the introduction, the average number $p$ of nonzero entries per
column and the number $p_{d}$ of nonzero entries in the densest column.
To make the efficiency comparison as fair as possible,
we have written the Fortran codes of SPAI, RSAI($tol$) and PSAI($tol$).
After $M$ is constructed by one of them, we then apply it within Krylov solvers.
We shall demonstrate that RSAI($tol$) works well. In the meantime, we compare
RSAI($tol$) with SPAI and PSAI($tol$), illustrating that
RSAI($tol$) is at least competitive with SPAI
and can be considerably more efficient and effective than the latter
for some problems, and it is as comparably effective as PSAI($tol$) for
preconditioning $Ax=b$.

\begin{table}[!htb]
\centering
\footnotesize
\tabcolsep 4pt
\caption{\label{table-mtr}The description of test matrices,
where $s$ is the number of irregular columns, $p$ the average number of
nonzero entries per column, and $p_d$ the number of the
nonzero entries in the densest column. An irregular column means that the
number of its nonzero entries exceeds $10p$.}
\begin{tabular}{ccccccc}\hline
matrices&$n$&$nnz$&$s$&$p$&$p_{d}$&Description\\
\hline
fs\_541\_3&541&4282&1&8&538&2D/3D problem\\
orsirr\_1&1030&6858&0&7&13&Oil reservoir simulation\\
orsirr\_2&886&5970&0&7&14&Oil reservoir simulation\\
sherman1&1000&3750&0&4&7&Oil reservoir simulation\\
sherman2&1080&23094&0&21&34&Oil reservoir simulation\\
sherman5&3312&20793&0&6&17&Oil reservoir simulation\\
saylr4&3564&22316&0&6&7&3D reservoir simulation\\
%Chebyshev2&2053&18447&0&8&9&Structural problem\\
cavity11&2597&71601&0&27&62&Subsequent computational fluid dynamics problem\\
ex36&3079&53099&0&17&37&Computational fluid dynamics problem\\
e20r0100&4241&131556&0&31&62&Computational fluid dynamics problem\\
e30r0000&9661&306356&0&32&62&Computational fluid dynamics problem\\
e40r0100&17281&553562&0&32&62&Computational fluid dynamics problem\\
powersim&15838&64424&2&4&41&Power network problem\\
raefsky3&21200&1488768&0&70&80&computational fluid dynamics problem\\
%hcircuit&105676&513072&151&5&1399&Circuit simulation problem\\
scircuit&170998&958936&104&6&353&Circuit, many parasitics\\
\hline
\end{tabular}
\end{table}

We perform numerical experiments on an Intel Core 2
Quad CPU E8400@ 3.00GHz with 2GB memory under
the Linux operating system. The computations of
constructing $M$ are done using Fortran 90 with the
machine precision $\epsilon_{\rm mach}=2.2\times10^{-16}$.
We take the initial sparsity pattern as that of
$I$ for SPAI and RSAI($tol$). We apply row
Dulmage-Mendelsohn permutations to the matrices
having zero diagonals so as to make their diagonals nonzero \cite{duff}.
The related Matlab commands are $j=\mathbf{dmperm}(A)$ and $A=A(j,:)$.
We applied $\mathbf{demperm}$ to cavity11, ex36, e20r0100, e30r0000 and
e40r0100. We use the $M$ generated by RSAI($tol$), SPAI and
PSAI($tol$) as right preconditioners, and use BiCGStab
as the Krylov solver, whose code is from Matlab 7.8.0.
The initial guess on the solution of $Ax=b$ is always $x_{0}=0$,
and the right-hand side $b$ is formed by choosing the solution
$x=(1,1,\ldots,1)^{T}$. The stopping criterion is
\begin{equation*}
  \frac{\|b-A\tilde{x}\|}{\|b\|}<10^{-8}, \tilde{x}=M\tilde{y},
\end{equation*}
where $\tilde{y}$ is the approximate solution obtained by BiCGStab applied to
the preconditioned linear system $AMy=b$.

In all the tables, $\varepsilon$, $l_{\max}$ and $c$ stand for
the accuracy requirements, the maximum loops
that RSAI($tol$) allows, and the numbers of the dominant
indices exploited by RSAI($tol$), respectively.
$spar=\frac{nnz(M)}{nnz(A)}$
denotes the sparsity of $M$ relative to $A$, $iter$ stands for the number of
iterations used by BiCGStab, $n_{c}$ is the number of columns of $M$ whose
residual norms do not drop below $\varepsilon$, and $ptime$ and
$stime$ denote the CPU timings (in seconds) of constructing $M$ and of
solving the preconditioned linear systems by BiCGStab, respectively.
$\dagger$ indicates that convergence is not attained within 1000 iterations,
$\ddag$ indicates that the Krylov solver stagnates and
$-$ means that we do not count CPU timings when BiCGStab fails to converge
within 1000 iterations.

\subsection{The effectiveness of the RSAI($tol$) algorithm} \label{subsec:5.1}
\begin{figure}[!htb]
\centering
\includegraphics{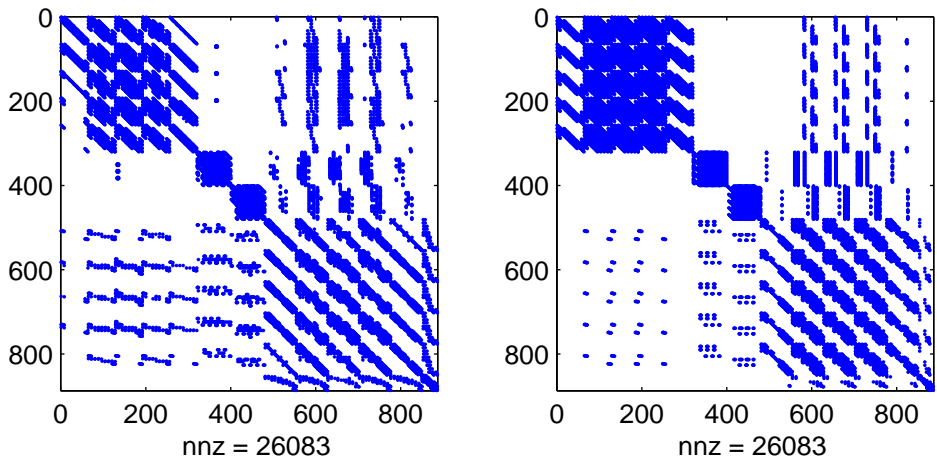}
\vspace{-80pt}
\caption{orsirr\_2: left: the pattern of $M$;
right: the sparsity pattern of $A^{-1}$.}\label{figure1}
\end{figure}

First of all, we illustrate that RSAI($tol$) can capture an effective approximate
sparsity pattern of $A^{-1}$ by taking the $886\times 886$ regular sparse matrix
orsirr\_2 as an example. We take the number $c=3$ of the dominant indices
exploited at each loop, the prescribed accuracy $\varepsilon=0.2$
and $l_{\max}=15$. We have found that all columns of $M$
satisfy the desired accuracy.
We use the Matlab code $inv(A)$ to compute $A^{-1}$
directly and then retain the $nnz(M)$ largest entries of $A^{-1}$
in magnitude. Figure \ref{figure1} depicts the patterns of $M$
and $A^{-1}$ that drops its small nonzero entries. We see that the
pattern of $M$ matches that of $A^{-1}$ quite well, demonstrating that
the RSAI($tol$) algorithm can indeed capture an effective approximate sparsity
pattern of $A^{-1}$. Importantly, it is clear from the figure that the sparsity
patterns of $M$ and $A^{-1}$ are irregular and the numbers of large entries
in the columns and rows of $A^{-1}$ vary greatly, although the matrix orsirr\_2
itself is regular sparse. As we have counted, the number of irregular columns
in $A^{-1}$ is 63, each of which has more than $10\frac{nnz(M)}{n}$ nonzero entries.
It means that about 8\% columns in $A^{-1}$ are irregular sparse.
This confirms the fact
addressed in \cite{Jia13b} that a good sparse approximate inverse
of a regular sparse matrix can be irregular sparse. Such fact also implies that
SPAI cannot guarantee to capture an effective approximate sparse pattern even
when $A$ is regular sparse and thus may be less effective for preconditioning
an regular sparse problem \eqref{eq1} since all the columns of $M$ constructed
by SPAI are sparse and each of them has at most $1+l_al_{\max}$ nonzero entries,
where $l_a$ is the number of most profitable indices added at each loop.

Now we take some difficult problems from Table \ref{table-mtr}
which make BiCGStab without preconditioning fail to converge within
1000 iterations. We precondition these problems by RSAI($tol$).
The preconditioners are all computed by setting
$\varepsilon=0.4$, $c=3$ and $l_{\max}=10$. Table \ref{table-di} lists the
results obtained.

\begin{table}[!htb]
\centering
\small
\caption{\label{table-di}The RSAI($tol$) algorithm for difficult problems}
\begin{tabular}{cccccc}\hline
matrices&$spar$&$ptime$&$n_{c}$&$iter$&$stime$\\
\hline
orsirr\_1&2.14&0.18&0&29&0.01\\
orsirr\_2&2.14&0.14&0&28&0.01\\
sherman2&2.76&1.59&0&6&0.01\\
sherman5&1.15&0.32&0&38&0.02\\
saylr4&1.22&0.92&0&672&0.49\\
ex36&2.40&4.82&10&90&0.11\\
e20r0100&2.77&20.42&148&148&0.45\\
raefsky3&1.46&104.70&0&206&4.78\\
\hline
\end{tabular}
\end{table}

We observe that the
RSAI($tol$) algorithm is effective to precondition these linear systems
and accelerates the convergence of BiCGStab dramatically in all cases,
as $iter$ indicates.
At the same time, we also find that the cost of constructing $M$ is
dominant. This is similar to SPAI and PSAI as well as
all the other non-factorized and factorized sparse approximate inverse
preconditioning procedures \cite{benzi02,benzi99}.
For all the problems, we see that the $n_{c}$
are equal to zero except for ex36 and e20r0100, for which $n_c$
is very small relative to $n$.
This means that for given parameters the RSAI($tol$) algorithm generally
construct effective preconditioners, as also reflected by
the $iter$, which are much smaller than 1000.
Therefore, we conclude from Table \ref{table-di} that the RSAI($tol$) algorithm
is generally effective for preconditioning \eqref{eq1}.

Next we vary the stopping criterion $\varepsilon$ to see how it
affects the sparsity of $M$ and convergence of BiCGStab. We still set
$c=3$ and $l_{\max}=10$. Table \ref{table-ep} reports the results
obtained for sherman1, saylr4 and orsirr\_2. As expected, BiCGStab used
fewer iterations in all cases and $M$ becomes denser as
$\varepsilon$ decreases. However, we find that, for the three problems,
though a smaller $\varepsilon$ makes BiCGStab converge faster,
it is more costly to construct a denser $M$. Compared with the cost
of constructing $M$, the cost of applying $M$ and solving the preconditioned
systems by BiCGStab is negligible for general problems, provided that
$M$ is an effective preconditioner. The choice
of $\varepsilon=0.4$ is the best as far as the total costs are concerned.
The table also implies that $l_{\max}=10$ is conservative
for the four given $\varepsilon$ since actual loops used does not achieve it
and whether or RSAI($tol$) terminated is up to $\varepsilon$ for the test
problems.

\begin{table}[!htb]
\centering
\small
\caption{\label{table-ep}The results for sherman1, saylr4
and orsirr\_2 with different $\varepsilon$}
\begin{tabular}{ccccc}\hline
\multicolumn{5}{c}{sherman1: $c=3$ and $l_{\max}=10$}\\
\hline
$\varepsilon$&$spar$&$ptime$&$iter$&$stime$\\
\hline
0.4&2.04&0.06&29&0.01\\
0.3&2.38&0.07&29&0.01\\
0.2&3.96&0.18&18&0.01\\
0.1&11.46&1.83&10&0.01\\
\hline

\multicolumn{5}{c}{saylr4: $c=3$ and $l_{\max}=10$}\\
\hline
$\varepsilon$&$spar$&$ptime$&$iter$&$stime$\\
\hline
0.4&1.22&0.92&672&0.49\\
0.3&1.95&1.90&267&0.21\\
0.2&3.08&3.84&85&0.08\\
0.1&6.79&12.50&45&0.05\\
\hline

\multicolumn{5}{c}{orsirr\_2: $c=3$ and $l_{\max}=10$}\\
\hline
$\varepsilon$&$spar$&$ptime$&$iter$&$stime$\\
\hline
0.4&2.14&0.14&28&0.01\\
0.3&2.75&0.21&23&0.01\\
0.2&4.36&0.52&16&0.01\\
0.1&8.23&2.35&11&0.01\\
\hline
\end{tabular}
\end{table}

Finally, we vary $c$ to investigate how it affects the sparsity of $M$ and
the convergence of BiCGStab. We set $\varepsilon=0.4$ and $l_{\max}=20$, and
take sherman2, ex36 and raefsky3 as examples.
Table \ref{table-mn} lists the results. We find that the $M$ gradually
become denser but the CPU time of constructing
them does not necessarily become more as $c$ increases. This should be expected
since at each loop the main cost of RSAI($tol$) is the solutions of LS problems
other than the ordering of entries of the current $r_k$ and
picking up indices.
By comparison, as far as the overall performance, measured by $ptime$, $n_c$ and
$iter$, is concerned, we find that  $c=3$ may be a very best choice.

\begin{table}[!htb]
\centering
\small
\caption{\label{table-mn} The results for sherman2, ex36 and raefsky3
with different $c$}
\begin{tabular}{cccccc}\hline
\multicolumn{6}{c}{sherman2: $\varepsilon=0.4$ and $l_{\max}=20$}\\
\hline
$c$&$spar$&$ptime$&$n_c$&$iter$&$stime$\\
\hline
1&1.99&1.70&2&8&0.01\\
2&2.13&1.55&0&7&0.01\\
3&2.76&1.59&0&6&0.01\\
4&3.22&1.70&0&7&0.01\\
5&3.55&1.94&0&7&0.01\\
\hline

\multicolumn{6}{c}{ex36: $\varepsilon=0.4$ and $l_{\max}=20$}\\
\hline
$c$&$spar$&$ptime$&$n_c$&$iter$&$stime$\\
\hline
1&1.83&7.12&59&103&0.11\\
2&2.22&5.75&0&85&0.11\\
3&2.40&4.86&0&91&0.11\\
4&2.67&4.94&0&83&0.11\\
5&2.93&5.00&0&86&0.11\\
\hline

\multicolumn{6}{c}{raefsky3: $\varepsilon=0.4$ and $l_{\max}=20$}\\
\hline
$c$&$spar$&$ptime$&$n_c$&$iter$&$stime$\\
\hline
1&1.17&120.49&0&$\dagger$&$-$\\
2&1.31&111.59&0&206&4.20\\
3&1.46&104.70&0&206&4.78\\
4&1.59&120.85&0&91&2.03\\
5&1.74&125.84&0&62&1.46\\
\hline
\end{tabular}
\end{table}

\subsection{The RSAI($tol$) algorithm versus the SPAI algorithm}
\label{subsec:5.2}

In this subsection, we compare the performance of RSAI($tol$) and
SPAI. For RSAI($tol$) we take $\varepsilon=0.3$,
the number $c=3$ of the dominant indices exploited at each loop
and $l_{\max}=10$. For SPAI we take the same $\varepsilon$ and
add three most profitable indices to the pattern of
$m_k$ at each loop for $k=1,2,\ldots,n$.
We set $l_{\max}=\lfloor10\times nnz(A)/(3\times n)\rfloor$ to control
$nnz(M)/nnz(A)\leq 5$ for the SPAI algorithm, where
$\lfloor x \rfloor$ is the Gaussian function.
Table \ref{table-rvs} presents the results.

\begin{sidewaystable}
%\begin{table}[!htb]
\centering
\caption{\label{table-rvs}The RSAI($tol$) algorithm versus
the SPAI algorithm}
\begin{tabular}{ccccccccccccccc}\hline
&&\multicolumn{6}{c}{The RSAI($tol$) algorithm}&&\multicolumn{6}{c}
{The SPAI algorithm}\\
\cline{3-8}\cline{10-15}
&&\multicolumn{3}{c}{Constructing $M$}&&\multicolumn{2}{c}{BiCGStab}&&
\multicolumn{3}{c}{Constructing $M$}&&\multicolumn{2}{c}{BiCGStab}\\
\cline{3-5}\cline{7-8}\cline{10-12}\cline{14-15}
matrices&&$spar$&$ptime$&$n_{c}$&&$iter$&$stime$&&$spar$&$ptime$&$n_{c}$
&&$iter$&$stime$\\
\hline
fs\_541\_3&&1.52&0.30&1&&16&0.01&&1.45&0.49&6&&18&0.01\\
orsirr\_1&&2.67&0.23&0&&24&0.01&&1.58&0.38&2&&29&0.01\\
orsirr\_2&&2.75&0.21&0&&23&0.01&&1.63&0.36&2&&26&0.01\\
sherman1&&2.38&0.07&0&&29&0.01&&1.72&0.12&7&&31&0.01\\
sherman2&&2.97&2.19&1&&5&0.01&&2.69&66.7&102&&$\dagger$&$-$\\
sherman5&&1.65&0.50&0&&30&0.02&&1.18&2.00&3&&34&0.02\\
saylr4&&1.95&1.90&0&&267&0.21&&1.97&2.09&15&&271&0.20\\
cavity11&&3.13&45.9&243&&172&0.28&&2.79&635.8&534&&215&0.34\\
ex36&&2.58&6.45&160&&77&0.10&&1.63&28.01&191&&93&0.10\\
e20r0100&&2.90&26.21&413&&149&0.47&&2.98&980.3&895&&212&0.68\\
e30r0000&&2.79&76.8&916&&168&1.28&&3.09&2727&2606&&211&1.69\\
e40r0100&&2.90&138.7&1664&&425&7.99&&3.06&5181&4649&&477&6.70\\
powersim&&4.92&24.2&1277&&839&3.02&&2.63&21.7&1414&&$\dagger$&$-$\\
raefsky3&&2.31&355&0&&75&2.12&&1.10&10655&94&&$\dagger$&$-$\\
scircuit&&2.82&1283&4&&849&46.41&&2.19&21134&7535&&$\dagger$&$-$\\
\hline
\end{tabular}
%\end{table}
\end{sidewaystable}

From the table, we see that the $M$ constructed by the RSAI($tol$) algorithm
are (almost) equally sparse as the counterparts by the SPAI algorithm
except powersim. However, the RSAI($tol$) algorithm uses substantially less
CPU time and more efficient than the SPAI algorithm except powersim, 
especially when a given $A$ is irregular sparse. For a dense column
of $A$, the cardinalities of the corresponding $\hat{\mathcal{J}}_k$ used by
SPAI are big during the loops, causing that it is time consuming to determine
the most profitable indices added to the patterns of $m_k$ at each loop.
In contrast, RSAI($tol$) overcomes this shortcoming by quickly finding new
indices added at each loop and is considerably more efficient. For sherman2,
raefsky3 and some others, we also find that the RSAI($tol$) algorithm
spends much less CPU timings than the SPAI algorithm when the given problem
is relatively dense, i.e., the average number $d$ of nonzero entries per column
is relatively large. This is because SPAI spent much more time seeking the most 
profitable indices for $d$ bigger than RSAI($tol$).

For the preconditioning effectiveness, we observe that the $n_{c}$
in RSAI($tol$) are considerably smaller than those in SPAI in all cases,
especially when $A$ is irregular sparse. The reason is that, as we have
explained, RSAI($tol$) is more effective to capture a good
approximate sparsity pattern of $A^{-1}$ than SPAI. Our results 
illustrate that with the roughly same number of nonzero entries allowed, 
the $M$ by RSAI($tol$) are more effective than their counterparts by SPAI, 
that is, RSAI($tol$) indeed captures the positions of large entries more 
reliably than SPAI does, and the latter may retain positions of some 
small nonzero entries but misses some large entries.  As a result,
the $M$ generated by SPAI are less effective than those by RSAI($tol$).
This is also confirmed by the number $iter$ of iterations, where the
$iter$ preconditioned by RSAI($tol$) are considerably fewer than those by
SPAI for all the test problems. Particularly,
for sherman2, powersim, raefsky3 and scircuit, SPAI fails to
make BiCGStab converge within 1000 iterations,
while RSAI($tol$) works well, and this is even true for the regular sparse
matrix sherman2 and raefsky3. Actually, for powersim and scircuit,
BiCGStab preconditioned by SPAI reduces the relative
residual to the level of $10^{-3}$ after several iterations, and
afterwards it does not decrease further.

In view of the above, the RSAI($tol$) algorithm is at least competitive
with and can be considerably more efficient and
effective than the SPAI algorithm, especially for the problems where
$A$ is irregular sparse or has relatively more nonzero entries.

\subsection{The RSAI($tol$) algorithm versus the PSAI($tol$) algorithm}
\label{subsec:5.3}

Keeping in mind the results of Section \ref{subsec:5.2}, we now compare the
RSAI($tol$) algorithm with the PSAI($tol$) algorithm proposed in
\cite{Jia09} and improved in \cite{Jia13a}. We also take $\varepsilon=0.3$
and $l_{\max}=10$ for PSAI($tol$). Table \ref{table-rvp} reports the
results obtained by PSAI($tol$).

\begin{table}[!htb]
\centering
\small
\caption{\label{table-rvp}The results by the PSAI($tol$)
algorithm}
\begin{tabular}{cccccccc}\hline
&&\multicolumn{3}{c}{Preconditioning}&&\multicolumn{2}{c}{BiCGStab}\\
\cline{3-5}\cline{7-8}
matrices&&$spar$&$ptime$&$n_{c}$&&$iter$&$stime$\\
\hline
fs\_541\_3&&1.87&0.13&0&&5&0.01\\
orsirr\_1&&5.36&1.69&0&&25&0.01\\
orsirr\_2&&5.66&1.55&0&&23&0.01\\
sherman1&&2.89&0.12&0&&27&0.01\\
sherman2&&2.74&1.59&0&&4&0.01\\
sherman5&&1.57&0.46&0&&29&0.02\\
saylr4&&1.86&12.26&0&&307&0.26\\
cavity11&&11.84&155.9&0&&55&0.35\\
ex36&&6.45&13.50&0&&55&0.12\\
e20r0100&&9.44&177.4&0&&91&0.79\\
e30r0000&&13.19&1670&29&&$\dagger$&$-$\\
e40r0100&&18.38&8779&494&&$\ddag$&$-$\\
powersim&&9.91&72.9&399&&52&0.41\\
raefsky3&&5.03&1281&0&&63&3.19\\
\hline
\end{tabular}
\end{table}

We observe that the $n_{c}$ by the PSAI($tol$) algorithm are no more than
those by the RSAI($tol$) algorithm for all the problems. Actually,
PSAI($tol$) obtains the $M$ with the desired accuracy satisfied for all
the problems except e30r0000, e40r0100 and powersim. This means that PSAI($tol$)
can construct effective preconditioners for most general sparse problems.
We do not list the test problem scircuit
in the table because PSAI($tol$) was found out of memory, which is due to
the occurrence of some large sized LS problems during the process.
We refer the reader to \cite{Jia13b,Jia09} for explanations
on some typical features of PSAI($tol$).
Compared with RSAI($tol$), it is seen that, except powersim, e30r0000 and
e40r0100, for the other test problems, there is no obvious winner
between it and PSAI($tol$) in term of $iter$.
Moreover, the setup time of $M$ by two algorithms are
also comparable for regular problems. For some problem which has relatively
more nonzero entries, such as raefsky3 and so on, the setup
time of $M$ by PSAI($tol$) is more than those by
RSAI($tol$). This is because we get a denser
$M$ by PSAI($tol$). However, it is more effective
than those by RSAI($tol$) in term of $n_{c}$.
Based on these observations, we conclude that the
RSAI($tol$) algorithm is as comparably effective as the
PSAI($tol$) algorithm.

\section{Conclusions}\label{sec:6}

SPAI may be costly to seek approximate sparsity patterns of $A^{-1}$
and ineffective for preconditioning a large sparse linear system,
which is especially true when $A$ is irregular sparse or is not very sparse. 
To this end, we have proposed a basic RSAI algorithm
that only uses the dominant information on residual and can
adaptively determine a good approximate sparsity pattern of $A^{-1}$.
We have derived an estimate for the number of
nonzero entries of $M$. In order to control the sparsity of $M$ and
improve efficiency, we have developed a practical RSAI($tol$) algorithm
with the robust adaptive dropping strategy
\cite{Jia13a} exploited. We have tested a number of real-world problems to
illustrate the efficiency and preconditioning effectiveness of RSAI($tol$). 
Meanwhile, we have numerically compared it with SPAI and PSAI($tol$), showing 
that RSAI($tol$) is at least competitive
with and can be substantially more efficient and effective than
SPAI, and it is also comparably effective as PSAI($tol$).

Some other research is significant. As we have seen,
RSAI($tol$) may be costly for $A$ row irregular sparse, which is also
the case for SPAI. In order to make RSAI($tol$) efficient 
for both column and row irregular sparse problems, we have exploited the 
idea proposed in \cite{Jia13b} to transform such a problem into certain 
regular ones, so that RSAI($tol$) can be much more
efficient to construct effective preconditioners \cite{jia15}. 
%It may be also worthwhile to nontrivially adapt RSAI($tol$) to
%large discrete ill-posed problems \cite{hansen98,hansen10}.
%\bigskip

%{\bf Acknowledgements}.

\small

\end{document}